\documentclass[12pt]{article}
\usepackage[centertags]{amsmath}
\usepackage{amsfonts}
\usepackage{amssymb}
\usepackage{amsthm}
\usepackage{newlfont}
\usepackage{graphics}
\usepackage{mathrsfs}
\usepackage{hyperref}

 \theoremstyle{definition}
 
 \theoremstyle{remark}

\newcommand{\vm}[1]{\mbox{\boldmath$#1$}}

\newcommand{\eps}{\varepsilon}
\newcommand{\sqep}{\sqrt{\varepsilon}}

\newcommand{\ds}{\displaystyle}
\newcommand{\beq}{\begin{eqnarray}}
\newcommand{\beqq}{\begin{eqnarray*}}
\newcommand{\eeq}{\end{eqnarray}}
\newcommand{\eeqq}{\end{eqnarray*}}
\newcommand{\x}{\mbox{\boldmath$x$}}
\newcommand{\y}{\mbox{\boldmath$y$}}
\newcommand{\z}{\mbox{\boldmath$z$}}
\newcommand{\w}{\mbox{\boldmath$w$}}

\newcommand{\e}{\mbox{\boldmath$e$}}

\newcommand{\pp}{\partial}

\begin{document}
\hsize=6 in  \textheight=9 in
\bibliographystyle{plain}
\baselineskip=24 pt
\title{ \textbf{ABOUT THE TRUE TYPE OF SMOOTHERS}}
\author {D. Ezri\thanks{Department of Electrical Engineering--Systems, Tel-Aviv University, Ramat-Aviv, Tel-Aviv 69978,
Israel. email: ezri@eng.tau.ac.il }   \and B. Z.
Bobrovsky\thanks{Department of Electrical Engineering--Systems,
Tel-Aviv University, Ramat-Aviv, Tel-Aviv 69978, Israel. email:
bobrov@eng.tau.ac.il} \and Z. Schuss\thanks{Department of
Mathematics, Tel-Aviv University, Ramat-Aviv, Tel-Aviv 69978,
Israel. email: schuss@post.tau.ac.il. The research of Z.S. was
partially supported by a research grant from TAU. } }

\maketitle

\begin{abstract}

We employ the variational formulation and the Euler-Lagrange
equations to study the steady-state error in linear non-causal
estimators (smoothers). We give a complete description of the
steady-state error for inputs that are polynomial in time. We show
that the steady-state error regime in a smoother is similar to that
in a filter of double the type. This means that the steady-state
error in the optimal smoother is significantly smaller than that in
the Kalman filter. The results reveal a significant advantage of
smoothing over filtering with respect to robustness to model
uncertainty.
\end{abstract}


\section{Introduction}\label{Introduction}
The appearance of steady-state errors is  endemic in linear and
nonlinear tracking systems and was extensively studied in control
theory \cite{Bishop}. An important result from linear deterministic
control theory is that the tracking properties of a closed-loop
control system are determined by its \emph{type}. The type of the
system is defined as the number of pure integrators in the transfer
function of the system within the closed loop. A type-$n$ system can
track a polynomial input of the form $a\ds\frac{t^{(p-1)}}{(p-1)!}$
without any steady-state error for $p \leq n$, with constant
steady-state error for $p=n+1$, and cannot track the input signal
for $p > n+1$. This result is easily obtained using the final value
theorem \cite{Bishop} from Laplace transform theory.

Another result from estimation theory is that a Kalman filter
\cite{Kalman_and_Bucy}, designed for the $n$-th order model
\begin{eqnarray}\label{NorderModelChap2}
\begin{array}{lll}
\ds\frac{d^n x}{dt^n}&=&x^{(n)}=\sigma\, \dot{w} \\&&\\
y&=&x+\rho\, \dot{v},
\end{array}
\end{eqnarray}
where $w(t)$ and $v(t)$ are independent Brownian motions, results in
a type-$n$ closed-loop tracking system (as demonstrated below).
Combining the above results  and using the linearity of the Kalman
filter, implies that a Kalman filter designed for the $n$-th order
model (\ref{NorderModelChap2}) tracks the signal
\begin{eqnarray}\label{XtildeNormalA}
\tilde{x}^{(n)}=\sigma\, \dot{w}+ a\frac{t^m}{m!}
\end{eqnarray}
with constant steady state error for $m=0$, and cannot track the
input signal for $m>0$.

Although a steady-state error is a fundamental concept in tracking
and control theory, the appearance of steady-state errors in
non-causal estimators (smoothers) has never been addressed, despite
the extensive study of linear and nonlinear smoothers  in the
literature
\cite{Rauch1,Rauch2,Rauch3,Lee,Weaver,bryson,BellmanKalaba,Sage,SageEwing,kailath22,Leondes,Anderson}.
The purpose of this paper is to fill the gap in the theory and to
give a complete description of the steady-state error regime in
linear smoothers. Specifically, we compute the steady-state error in
the linear minimum mean square error (MMSE) smoother, designed for
the model (\ref{NorderModelChap2}). We show that the steady-state
error in the smoother is similar to that in the corresponding filter
of double the type. Questions of convergence and stability of fixed
interval smoothers were considered recently in \cite{Einicke}.

\section{Mathematical preliminaries}
We begin with the general linear signal and measurement model
\begin{eqnarray}\label{LinearModelIntro}
\begin{array}{lll}
\dot{\x}&=&\vm{A}\x+\vm{G}\dot{\w} \\ && \\
\y & = & \vm{H}\x+\vm{B}\dot{\vm{v}},
\end{array}
\end{eqnarray}
where $\w(t), \ \vm{v}(t)$ are orthogonal vectors of independent
Brownian motions. The minimum mean square estimation error (MMSEE)
or maximum a-posteriori (MAP) estimator $\hat{\x}(\cdot)$ of
$\x(\cdot)$, conditioned on the measurement $\y(\cdot)$, is the
minimizer over all square integrable functions $\vm{u}$ of the
integral quadratic cost functional \cite{bryson,BellmanKalaba}
\begin{eqnarray}\label{LQGFunctional}
\int_0^T\left[
(\y-\vm{Hz})^T(\vm{BB}^T)^{-1}(\y-\vm{Hz})+\vm{u}^T\vm{u}\right]\,dt,
\end{eqnarray}
subject to the equality constraint
\begin{eqnarray}\label{LinearEqaulityConstraint}
\dot{\z}=\vm{Az}+\vm{Gu},
\end{eqnarray}
where $[0,T]$ is the observation interval.  Note that the first term
in (\ref{LQGFunctional}) represents the energy of the measurements
noise associated with the test function $\vm{z}(\cdot)$, while the
second term represents that of the driving noise associated with
$\vm{z}( \cdot)$. Note further that the integral in
(\ref{LQGFunctional}) contains the white noises $\dot{\w}(t),\
\dot{\vm{v}}(t)$, which are not square integrable. To remedy this
problem, we proceed in the standard way \cite{stroock} by beginning
with a model in which the white noises $\dot{\w}(t),\
\dot{\vm{v}}(t)$ are replaced with square integrable wide band
noises, and at the appropriate stage of the analysis, we take the
white noise limit.

The Euler-Lagrange (EL) equations \cite{Kirk} for the minimizer
$\hat{\x}=\hat{\x}(t\,|\,T)$ of the problem (\ref{LQGFunctional}),
(\ref{LinearEqaulityConstraint}) are
\begin{eqnarray}\label{GeneralEL}
\begin{array}{lll}
\ds\frac{\pp\hat{\x}(t\,|\,T)}{\pp t}&=&\vm{A} \hat{\x}(t\,|\,T)
-\ds\frac{1}{2}\vm{GG}^T\vm{\lambda} \\
&&\\
\dot{\vm{\lambda}}&=&-2\vm{H}^T(\vm{BB}^T)^{-1}\vm{H}\hat{\x}(t\,|\,T)-\vm{A}^T\vm{\lambda}
+2\vm{H}^T(\vm{BB}^T)^{-1}\y(t),
\end{array}
\end{eqnarray}
with the boundary conditions
\begin{eqnarray}
\hat{\x}(0\,|\,T)=\x_0, \ \ \ \ \ \vm{\lambda}(T)=0.
\end{eqnarray}
Thus, the estimation problem is reduced to a  linear two-point
boundary-value problem, which is solved in closed form by the sweep
method. Specifically, first the signal is filtered causally by the
Kalman filter $\hat{\x}(t)$, satisfying
\begin{eqnarray}\label{KalmanEquation}
\dot{\hat{\x}}(t)&=&\vm{A}\hat{\x}(t)+\vm{P}\vm{H}^T(\vm{BB}^T)^{-1}[\y(t)-\vm{H}\hat{\x}(t)],
\end{eqnarray}
and $\vm{P}(t)$ is the covariance matrix, satisfying the Riccati
equation
\begin{eqnarray}\label{RiccatiForKalman}
\dot{\vm{P}}=\vm{AP}+\vm{PA}^T+\vm{GG}^T-\vm{PH}^T(\vm{BB}^T)^{-1}\vm{HP}.
\end{eqnarray}
Then the smoother $\hat{\x}(t\,|\,T)$ is the backward sweep of the
Kalman filter,
\begin{eqnarray}\label{ForwardSolutionForSmoother}
\ds\frac{\pp\hat{\x}(t\,|\,T)}{\pp
t}=\vm{A}\hat{\x}(t\,|\,T)+\vm{GG}^T\vm{P}^{-1}[\hat{\x}(t\,|\,T)-\hat{\x}(t)],
\end{eqnarray}
with the boundary condition
\begin{equation}
\hat{\x}(T\,|\,T)=\hat{\x}(T).
\end{equation}

We continue with the steady-state error regime in the Kalman filter
designed for the $n$-th order model (\ref{NorderModelChap2}).
Rewriting (\ref{NorderModelChap2}) in vector matrix notation leads
to
\begin{eqnarray}\label{VMNorderModelChap2}
\begin{array}{lll}
\dot{\x}&=&\vm{A}\x+\sigma\vm{b}\dot{w} \\&&\\
y&=&\vm{h}^T \x+\rho\,\dot{v},
\end{array}
\end{eqnarray}
where $\x=[x, x_2, \ldots, x_n]^T$, and
\begin{eqnarray*}
\vm{A}=\left [ \begin{array}{cccccc}
  0 & 1 & 0 & 0 & \ldots & 0 \\
  0 & 0 & 1 & 0 & \ldots & 0 \\
  0 & 0 & 0 & 1 & \ldots & 0 \\
    &   & \vdots &   &   &   \\
  0 & 0 & 0 & 0 & \ldots & 0
\end{array}  \right];\quad \vm{b}=\left [ \begin{array}{c}
  0 \\
  0 \\
  \vdots \\
  0 \\
  1
\end{array}  \right ];\quad \vm{h}=\left [ \begin{array}{c}
  1 \\
  0 \\
  \vdots \\
  0 \\
  0
\end{array}  \right ].
\end{eqnarray*}
The Kalman filter estimator $\hat{\x}(t)=[\hat{x}(t),
\hat{x}_2(t),\ldots,\hat{x}_n(t)]^T$ for the model
(\ref{VMNorderModelChap2}) satisfies (see equation
(\ref{KalmanEquation})
\begin{eqnarray}\label{VMKalmanChap2}
\dot{\hat{\x}}(t)=\vm{A}\hat{\x}(t)+\vm{k}[y(t)-\hat{x}(t)],
\end{eqnarray}
where
$\vm{k}(t)=[k_1(t),\ldots,k_n(t)]^T=\ds\frac{1}{\rho^2}\vm{P}(t)\vm{h}$.
The equation of the Kalman filter (\ref{VMKalmanChap2}) may be
rewritten in the form
\begin{eqnarray}
\begin{array}{lll}
\dot{\hat{x}}(t)&=&\hat{x}_2(t)+k_1[y(t)-\hat{x}(t)] \\
&\vdots& \\
\dot{\hat{x}}_{n-1}(t)&=&\hat{x}_n(t)+k_{n-1}[y-\hat{x}(t)] \\&&\\
\dot{\hat{x}}_{n}(t)&=&k_{n}[y-\hat{x}(t)],
\end{array}
\end{eqnarray}
which form a type-$n$ closed-loop tacking system with input $y(t)$,
output $\hat{x}(t)$ and open-loop transfer function
\begin{eqnarray}
G(S)=\frac{k_1 S^{n-1}+k_2 S^{n-2}+\ldots+k_n}{S^n}.
\end{eqnarray}
Thus, the Kalman filter designed for the $n$-th order model
(\ref{NorderModelChap2}) tracks the signal
\begin{eqnarray}\label{XtildeNormalA}
\tilde{x}^{(n)}=\sigma\, \dot{w}+ a\frac{t^m}{m!}
\end{eqnarray}
with constant steady state error for $m=0$, and cannot track the
input signal for $m>0$. The the terminology \emph{offset} case $m=0$
of filtering theory is kept also for smoothing.

\section{Steady-state errors in smoothers}

Having described the steady-state error phenomenon in the Kalman
filters, we continue with the evaluation of the steady-state errors
in smoothers. We assume that a linear smoother is designed for the
$n$-th order signal model (\ref{NorderModelChap2}), for which the
Kalman filter was designed in the previous section. In order to
examine the tracking properties of the smoother, we assume that the
incoming signal is augmented with a $(n+m)-$th order polynomial in
$t$
\begin{eqnarray}\label{LinearSignalModelWithTension}
\begin{array}{lll}
  \tilde{x}^{(n)} & = & \ds\frac{a t^m}{m!} +\sigma \,\dot{w} \\
  &&\\
  y & = & \tilde{x}+\rho\,\dot{v}.
\end{array}
\end{eqnarray}
For the sake of simplicity, we assume hereafter that
$\sigma=\rho=\sqep$ and $\tilde{x}(0)=\hat{x}(0)=0$. The smoother is
a linear system, so the response of the smoother to the measurement
signal $y(t)$ can be separated into a deterministic term, depending
on the drift $\ds\frac{a t^m}{m!} $, and a stochastic term depending
on $w(t),\,v(t)$. Because the steady-state error regime is
determined by the deterministic term, we consider the noiseless
version of (\ref{LinearSignalModelWithTension})
\begin{eqnarray}\label{NoiselessLinearSignalModelWithTension}
\begin{array}{lll}
  \tilde{x}^{(n)} & = & \ds\frac{a t^m}{m!}  \\
  &&\\
  y & = & \tilde{x}.
\end{array}
\end{eqnarray}
Using vector matrix notation, the smoother is designed for the model
(\ref{VMNorderModelChap2}) and the equations for the smoothed
estimate $\hat{\x}(t\,|\,T)$ (\ref{ForwardSolutionForSmoother}),
(\ref{KalmanEquation}) are
\begin{eqnarray}\label{LinearForwardSweep}
\begin{array}{llll}
&&&\\
 \ds\frac{\pp{\hat{\x}}(t\,|\,T)}{\pp t}&=&\vm{A}\hat{\x}(t\,|\,T)+\ds
\frac{\vm{Ph}}{\eps}\left
[\y(t)-\vm{h}^T\hat{\x}(t\,|\,T)\right ], & \hat{\x}(0)=0  \\
&&&\\
 \dot{\hat{\x}}(t\,|\,T))&=&\vm{A}\hat{\x}(t\,|\,T)+\vm{bb}^T\vm{P}^{-1}\left
[\hat{\x}(t\,|\,T)-\hat{\x}(t)\right ], & \hat{\x}(T|T)=\hat{\x}(T).
\end{array}
\end{eqnarray}

Inspecting (\ref{LinearForwardSweep}), we observe that the Kalman
filter forms a feedback system with a type-$n$ transfer function.
Thus, the Kalman filter is capable of tracking a $(p < n)$-th order
polynomial in $t$ without a steady-state error, and a $(p =n)$-th
order polynomial in $t$ with a constant steady-state error. The
backward equation for $\hat{\x}(t\,|\,T)$ forms a feedback system
tracking the output of the Kalman filter $\hat{\x}(t)$ in reverse
time and the transfer function is again of type-$n$, with tracking
properties similar to these of $\hat{\x}(t)$. An important feature
of this forward-backward tracking system is that the tracking errors
in the backward and forward equations have reverse signs. This
qualitative observation may indicate that the steady-state error
regime in the smoother is superior to that in the causal filter.

With the above observations in mind,  we turn to the quantitative
analysis of the steady-state error in the smoother. We define the
estimation error $\vm{e}(t\,|\,T)=\hat{\x}(t\,|\,T)-\tilde{\x}(t)$
and use the expression (\ref{NoiselessLinearSignalModelWithTension})
for the incoming signal, which we rewrite as
\begin{eqnarray}
\begin{array}{lll}
&& \\
\dot{\tilde{\x}}(t)&=&\vm{A\tilde x}(t)+\ds\frac{at^m}{m!} \vm{b}\\
&& \\
y(t)&=&\vm{h}^T \tilde{\x}(t),  \\
&&
\end{array}
\end{eqnarray}
where $\vm{A},\,\vm{b}$ and $\vm{h}$ are given in
 (\ref{VMNorderModelChap2}).  Using the EL equations (\ref{GeneralEL}), the equations
for the estimation error are
\begin{eqnarray}\label{ErrorEquationForSteadyStateErrorMatrix}
\begin{array}{lllll}
\ds\frac{\pp{\vm{\e}}(t\,|\,T)}{\pp t}&=& \vm{Ae}(t\,|\,T)-\ds
\frac{1}{2}\vm{bb}^T
\vm{\lambda}(t)-\frac{at^m}{m!},&& \vm{e}(0\,|\,T)=0  \\
&&& \\  \dot{\vm{\lambda}}(t)&=& -2 \vm{hh}^T \vm{e}(t\,|\,T)
-\vm{A}^T\vm{\lambda}(t), && \vm{\lambda}(T)=0.
\end{array}
\end{eqnarray}
Setting $\vm{e}(t\,|\,T)=[e_1, \ldots, e_n]^T,\
\vm{\lambda}(t)=[\lambda_1, \ldots, \lambda_n]^T$, the error
equations (\ref{ErrorEquationForSteadyStateErrorMatrix})  become
\begin{eqnarray}\label{EulerLagrangeForSteadyStateErrorSeparated}
\begin{array}{lllll}
&&&&\\
\dot{e}_k&=&e_{k+1}, &&\ e_k(0\,|\,T)=0 \quad\mbox{for all}\quad 1 \leq k \leq n-1 \\
&&&&\\
\dot{e}_n & = &  -\ds \frac{1}{2}\lambda_n-\ds \frac{at^m}{m!}, &&
e_n(0\,|\,T)=0 \\
&&&& \\
\dot{\lambda}_1 & = & -2 e_1, &&  \lambda_1(T)=0 \\
&&&& \\
\dot{\lambda}_k & = & -\lambda_{k-1}, && \lambda_k(T)=0 \quad\mbox{for all}\quad 2 \leq k \leq n. \\
&&&&
\end{array}
\end{eqnarray}
We seek a steady-state solution, such that $\dot{e}_k =0
\quad\mbox{for all}\quad 1 \leq k \leq n$. This implies
\begin{eqnarray}
\dot{e}_n & = & -\frac{1}{2}\lambda_n-\frac{at^m}{m!} = 0,
\end{eqnarray}
hence
 \[ \lambda_n=-\frac{2at^m}{m!}.\]
Using (\ref{EulerLagrangeForSteadyStateErrorSeparated}), we have for
all  $0 \leq k \leq n-1$
\begin{eqnarray}
\dot{\lambda}_{n-k} = \left \{
\begin{array}{ccc}
  (-1)^{k+1} \ds\frac{2a\,t^{m-k-1}}{(m-k-1)!} && \mbox{for all} \quad k \leq m-1 \\
  & &\\
  0 && \mbox{for all} \quad k > m-1.
\end{array}  \right.
\end{eqnarray}
Using the fact that $\dot{\lambda}_1=-2\,e_1$, we obtain the
steady-state solution for $e_1$, whenever it exists, as
\begin{eqnarray}
e_1=\left \{
\begin{array}{cc}
  0 & \mbox{for all} \quad  m < n \\ & \\
  (-1)^{n}a & \mbox{for all} \quad m=n \\ & \\
  \mbox{no steady-state solution} & \mbox{for all} \quad m>n.
\end{array} \right.
\end{eqnarray}

\begin{figure}
\centering
\resizebox{!}{10cm}{\includegraphics{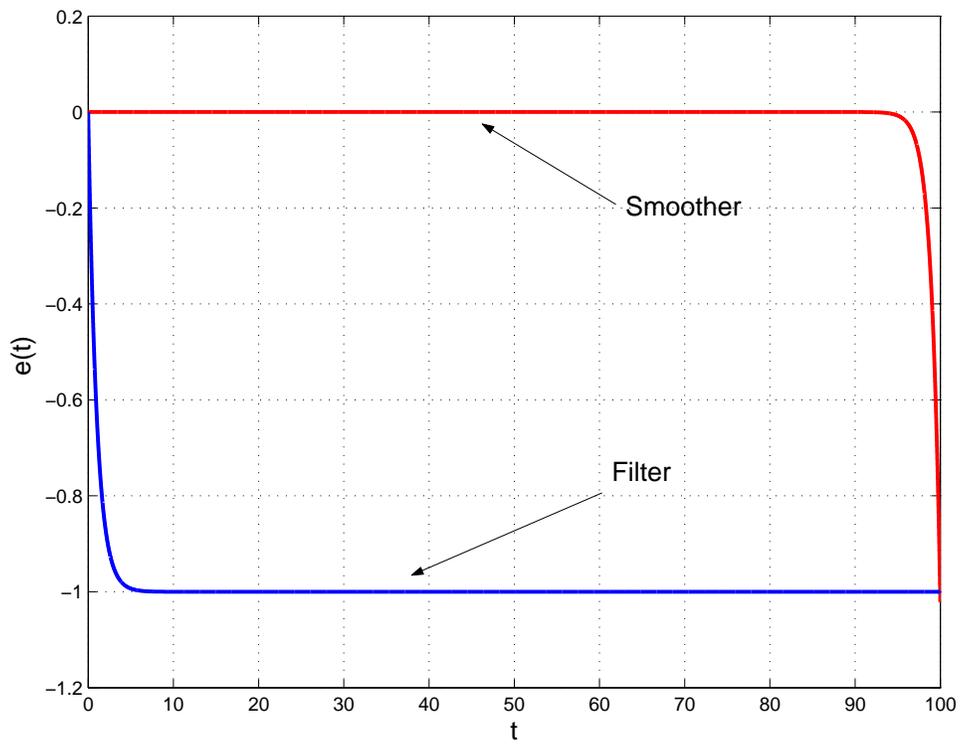}}
\caption{\label{Stead_State_Error_1offset}The steady-state error
in the first order smoother and filter with offset $(n=1,\ m=0)$.
No steady-state error develops in the smoother, while, the filter
develops constant steady-state error. }
\end{figure}

\begin{figure}
\centering
\resizebox{!}{10cm}{\includegraphics{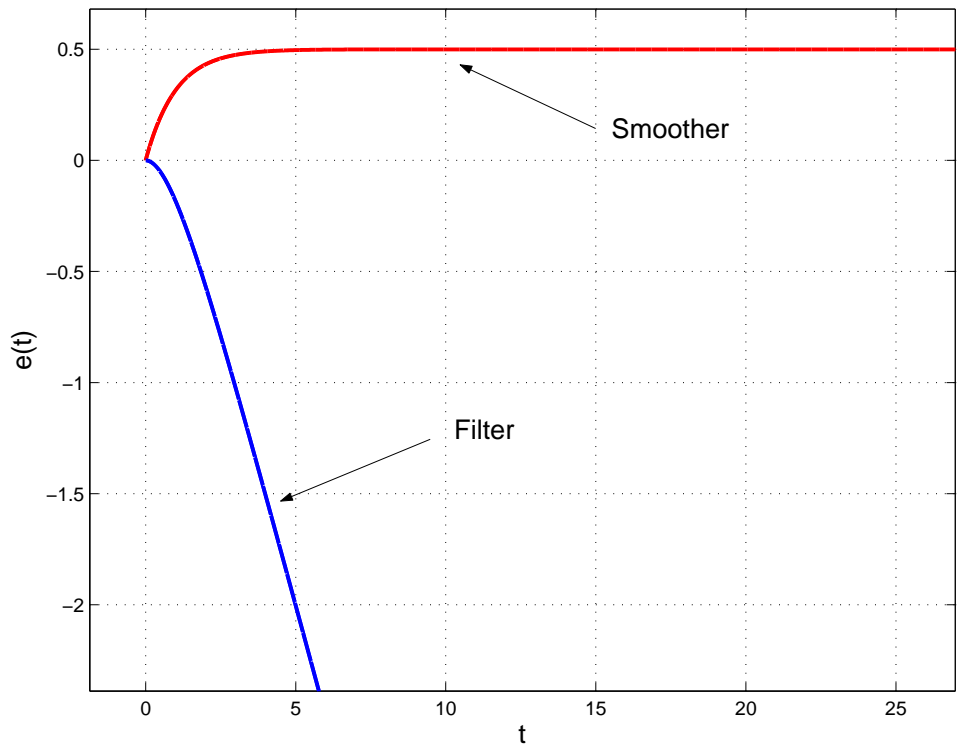}}
\caption{\label{Stead_State_Error_2offset}The steady-state error
of the first order smoother and filter with the polynomial input
$\ds\frac{at^2}{2}$ $(n=1,\ m=1)$. The smoother develops constant
steady-state error, while, the error in the filter is unbounded. }
\end{figure}

This result implies that a smoother designed for a $n$-th order
model is capable of tracking polynomial inputs of order up to
$2n-1$, that is, for $m=n-1$, without any steady-state error.
Moreover, for a $2n$-th order polynomial input (i.e., $m=n$) a
constant steady-state error of magnitude $|a|$ appears. We conclude
that the steady-state regime in a smoother designed for a $n$-th
order model is similar to a causal filter designed for a $2n$-th
order model, as illustrated in Figures
\ref{Stead_State_Error_1offset}, \ref{Stead_State_Error_2offset}.

\begin{figure}
\centering \resizebox{!}{10cm}{\includegraphics{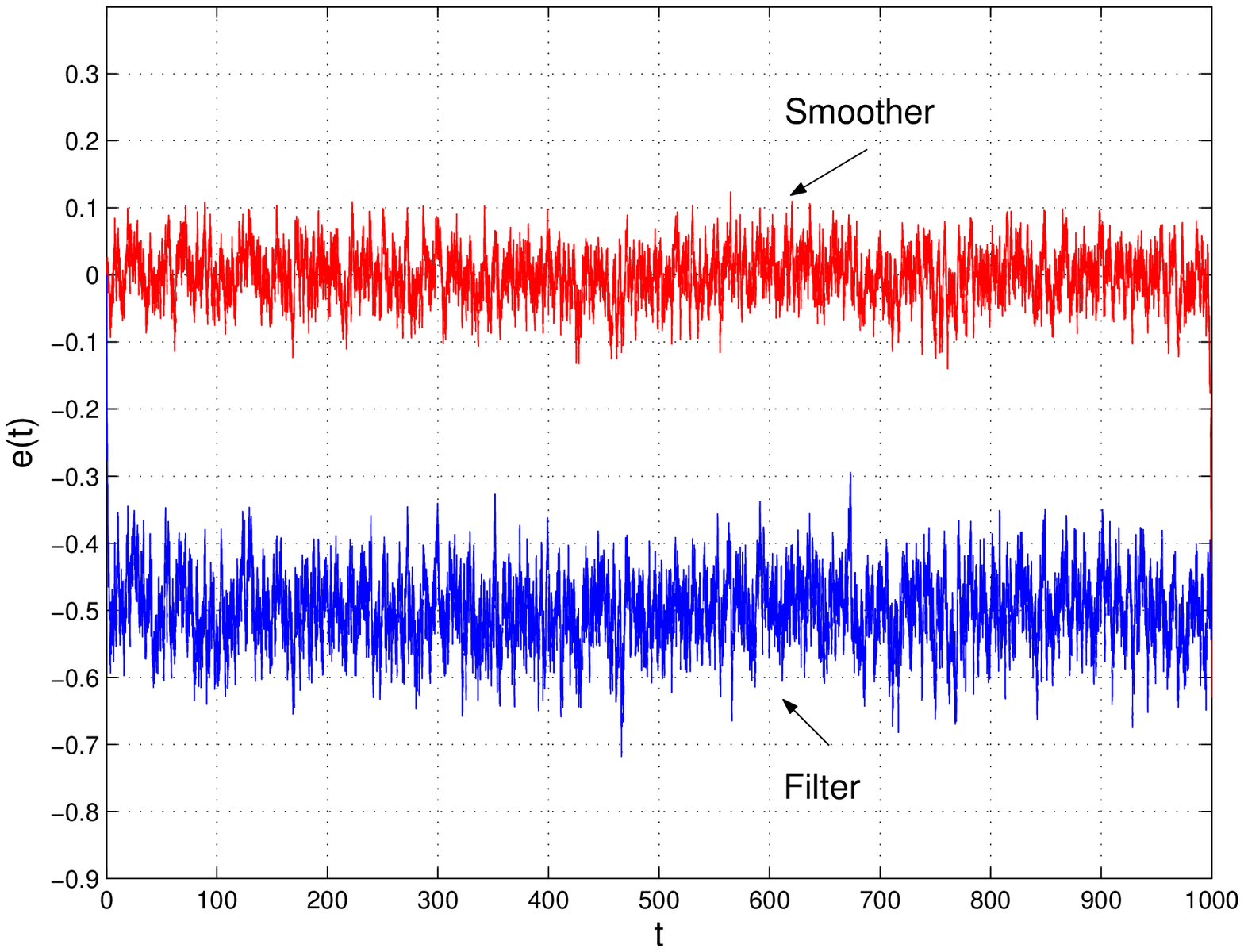}}
\caption{\label{ErrorWithOffset}The estimation error in the first
order smoother and filter in the offset case. No steady-state
error develops in the smoother and the error variance in the
smoother remains smaller.}
\end{figure}

The linearity of the smoother implies that  the smoothing error
variance induced by the Brownian motions $w(t), \, v(t)$ is not
affected by the drift term $\ds \frac{at^m}{m!}$ in the signal model
(\ref{LinearSignalModelWithTension}). Thus, the smoother exhibits
superior steady-state error regime without any increase in the error
variance, which remains smaller than that in the causal filter
\cite{kailath22}, as shown in Figure \ref{ErrorWithOffset}.

The obtained steady-state error regime corresponds to an infinite
interval $0<t<\infty$. Obviously, at the end of the interval $[0,T]$
the smoothing solution $\hat{\x}(T\,|\,T)$ is equal to the filtering
solution $\hat{\x}(T)$ (\ref{LinearForwardSweep}). Thus, in the
offset case ($m=0$) the filter develops a constant steady-state
error and the lag needed to eliminate the steady-state error by
smoothing is the time constant of the system. In contrast, when
$m>0$ and the causal filter error increases with time, the lag
necessary to decrease the error to the values we obtained is much
longer.


\section{Discussion and conclusions}
The results  reveal a significant robustness of the smoother to
signals different than those, for which the smoother is optimal.
In many applications in engineering practice, either there is no
exact information about the signal, or the signal changes its
characteristics with time. Thus, robustness of the system to the
signal model is desirable. The results further indicate that the
performance gap between smoothers and filters with respect to the
criterion of mean square error increases as the difference between
the incoming signal and the nominal signal increases.

The fact that a smoother, optimal for a $n$-th order signal,
exhibits a steady-state error regime of a type $2n$ causal filter,
without degradation in the error variance, is very appealing from
the engineering point of view. This is due to the facts that in many
applications lag estimation is possible and that linear smoothers
are rather easy to implement, for example, by using approximate
finite impulse response (FIR) non-causal filters. Utilizing
smoothers to exploit their superior tracking properties may free the
designer from the traditional trade-off that exists in causal
filtering between steady-state error regime and error variance.

Steady-state errors also appear in nonlinear estimators, as shown in
nonlinear filtering \cite{WeltiJsac,WeltiThird}. An extension of the
ideas demonstrated here to  nonlinear smoothing will be given in a
forthcoming paper.

\end{document}